\documentclass[]{article}
\usepackage{enumerate}
\usepackage{amsmath,amsthm}
\usepackage{amsfonts}
\usepackage{amssymb}

\begin{document}
\title{Recurrence and transience of near-critical multivariate growth models: criteria and examples}
\author{G\"otz Kersting\thanks{Institut f\"ur Mathematik, Goethe Universit\"at, Frankfurt am Main, Germany, kersting@math.uni-frankfurt.de}}
\date{\today}
\maketitle

\begin{abstract} 
We discuss complementary recurrence and transience criteria for stochastic processes $(X_n)_{n \ge 0}$ with values in the $d$-dimensional orthant $\mathbb R^d_+$ fulfilling
a non-linear stochastic equation of the form $X_{n+1}=MX_n+g(X_n)+ \xi_n$ with a primitive matrix $M$ and random noise $\xi_n$ and obeying a weak Markov property. As examples we discuss bisexual Galton-Watson processes and multivariate Galton-Watson processes, which both may be population size dependent.
\medskip



\end{abstract}

\section{Introduction}

We consider stochastic processes $(X_n)_{n\ge 0}$ taking values in the $d$-dimensional orthant $\mathbb R_+^d= \{ (x_1, \ldots,x_d)^T \in \mathbb R^d : x_i \ge 0 \}$ and adapted to some filtration $(\mathcal F_n)_{n \ge 0}$, which satisfy an equation of the form
\[ X_{n+1} = MX_n + g(X_n)+ \xi_n \ , \quad n \in \mathbb N_0\ ,   \]
with a $d\times d$ matrix $M$  having non-negative entries, with a measurable function $g: \mathbb R_+^d \to \mathbb R^d$, and with random fluctuations $\xi_n=(\xi_{n1}, \ldots, \xi_{nd})^T$ satisfying
\[ \mathbf E[\xi_n \mid \mathcal F_n]=0 \ \text{ a.s.} \]
One may view the process as a non-linear random pertubation of the linear dynamical system $x_{n+1}=Mx_n$, $n \in \mathbb N_0$. Here we require this system to be critical, which means that the Perron-Frobenius eigenvalue $\lambda_1$ of $M$ is equal to 1. We focus on the situation, when $M$ is a primitive matrix; then up to scaling there is a unique  left eigenvector $\ell=(\ell_1,\ldots, \ell_d)$ corresponding to $\lambda_1$ and its components are all strictly positive. As usual we let $\ell_1+ \cdots +\ell_d=1$.

Now the size of the random fluctuations will be determined by the conditional variance of $\ell \xi_n= \ell_1\xi_{n1} + \cdots + \ell_d\xi_{nd}$. More precisely we assume that
\[ \mathbf E[ (\ell \xi_n)^2 \mid \mathcal F_n] = \sigma^2(X_n) \  \text{ a.s.} \]
with some measurable function $\sigma^2: \mathbb R^d_+ \to \mathbb R_+$. 

We are aiming at general criteria for recurrence or transience of the process $(X_n)_{n \ge 0}$, that is, at criteria which allow to decide whether the event $\{\|X_n\| \to \infty$ for $n \to \infty\}$ is an event of zero probability or not, where $\| \cdot \|$ denotes an arbitrary norm on $\mathbb R^d$. Hereby, in talking about recurrence and transience, we have taken the liberty to adopt the  terminology from Markov chain theory. Certainly our processes obey only a relaxed form of the Markov property,  but examples typically are Markov chains.

Our theorems in the multivariate setting are complete generalizations of the known results in the univariate setting. Therefore it is appropriate to first reconsider the univariate case. This is done in section 2. In section 3 we apply these results to the population size dependent bisexual Galton-Watson process. The multivariate  case is then discussed in section 4. As an example the multivariate population size dependent Galton Watson process is treated in section 5. Proofs for the multivariate criteria are given elsewhere.

\section{The univariate case revisited}

In the 1-dimensional case our model equation simplifies to the difference equation
\[ X_{n+1}=X_n+g(X_n)+ \xi_n  \]
with some  function $g:\mathbb R_+ \to \mathbb R$. 
A number of examples can be put into this framework, among others e.g. population size dependent branching processes \cite{kle}, \cite{ku}, controlled branching processes \cite{go}, branching processes in random environment \cite{al}, or nonlinear stochastic trends \cite{gra}.

Our main condition is 
\begin{align} \label{A1} g^+(x) =o(x)  \ \text{ as } x \to \infty \ . \tag{A1} 
\end{align}
This  assumption of ``near-criticality''  simply says that $X_n$ is the dominating term within $X_n+g(X_n)$ such that supercritical growth is excluded.  

Also we assume the existence of some constants $c, \delta>0$ 
such that for $X_n \ge c$
\begin{align} \label{A2} \mathbf E[ |\xi_n|^p \mid \mathcal F_n] \le c \sigma^p(X_n)\ \text{ a.s. with } p=2+\delta \ . \tag{A2} \end{align}
 With these two conditions we have the following  criteria complementary to each other.

\paragraph{Theorem 1.} {\em Let \eqref{A1}, \eqref{A2} be fulfilled. Assume that there is an $\varepsilon >0$ such that
\begin{align}  xg(x) \le \frac {1-\varepsilon}2 \sigma^2(x) \label{recurr}
\end{align}
for $x$ sufficiently large. Then
\[ \mathbf P( X_n \to \infty \text{ for } n \to \infty)=0 \ . \]}

\noindent
The converse criterium requires some slight additional restrictions.

\paragraph{Theorem 2.}{\em  Let \eqref{A1}, \eqref{A2} be fulfilled. Also let $\sigma^2(x)$ be bounded away from zero on intervals $(u,v)$ with $0<u<v<\infty$, and let
\begin{align} \label{sigma} \sigma^2(x)= O(x^2 \log^{-2/\delta} x) \ \text{ for } x \to \infty  \end{align}
with $\delta$ as in \eqref{A2}. Assume that there is an $\varepsilon>0$ such that
\begin{align} \label{trans}  xg(x) \ge \frac {1+\varepsilon}2 \sigma^2(x) \end{align}
for $x$ sufficiently large. Then there is a number $m <\infty$ such that
\[ \mathbf P( \limsup_n X_n \le m \text{ or } \lim_n X_n = \infty)=1 \ . \]
If also for every $c>0$ there is a $n \in \mathbb N$ such that $\mathbf P(X_n>c) >0$, then
\[ \mathbf P(X_n \to \infty \text{ for } n \to \infty)>0 \ .\]}

These results are contained in \cite{ke}. From there the first theorem is taken literally, while the second one is a somewhat more general version of the corresponding Theorem 2 in \cite{ke}. There the condition $g(x)= O(x\log^{-2/\delta}x) $ is used, which  is stronger than our condition \eqref{sigma} in view of   \eqref{trans}. Thus our theorem offers a relaxation of conditions, which is more to the point and also useful in examples, while   the proof of the criterion remains practically unchanged (as one easily convinces oneself).  

The theorems can be understood as follows: Typically the long term behavior is either dominated by the ``drift term'' $g(X_n)$, or it is mainly controlled by the fluctuations $\xi_n$. There is only a small boundary region where both the drift term and the fluctuations have to be taken into account. It is there, where one would expect a particular rich and variable stochastic behavior. 

\paragraph{Remark 1.}  In order to get a better understanding of the main condition of both theorems it is instructive to rewrite the model equation in a multiplicative form as
\[ X_{n+1} =X_n(1+ h(X_n) + \zeta_n) \quad \text{with } h(x)=\frac {g(x)}x \ , \ \zeta_n= \frac{\xi_n}{X_n} \ . \]
Then
\[ \mathbf E[ \zeta_n^2 \mid \mathcal F_n] = \tau^2(X_n) \text{ a.s.}Ê \quad \text{with } \tau^2(x)= \frac {\sigma^2(x)}{x^2} \ . \]
Now the main requirements \eqref{recurr} and \eqref{trans} of the theorems read $h(x) \le \frac{1-\varepsilon}2 \tau^2(x)$ versus 
$h(x) \ge \frac{1+\varepsilon}2 \tau^2(x)$. In this formulation the drift is directly related to the variance of the fluctuations. 

\paragraph{Remark 2.} One might wonder, whether condition \eqref{sigma} can be  substantially  relaxed or even removed. Without compensation this is not possible, as can be seen from Example C in Section 3 of \cite{ke}.

\section{Example: The bisexual GW-process}

For the bisexual Galton-Watson process the $n$-th generation of some po\-pulation (with $n=0,1,\ldots$) consists of $F_n$ female and $M_n$ male individuals. They are assumed to form $L(F_n,M_n)$   different couples, with some given deterministic ``mating function'' $L(\cdot, \cdot)$, such that $L(0,\cdot)=L(\cdot,0)=0$. The $i$-th couple then has $\rho_{ni}$ female and $\tau_{ni}$ male offspring. Thus the population evolves according to the equations
\[ F_{n+1}= \sum_{i=1}^{L(F_n,M_n)} \rho_{ni} \ , \quad M_{n+1} = \sum_{i=1}^{L(F_n,M_n)} \tau_{ni} \ . \] 
Let $\mathcal F_n$ be the $\sigma$-field, generated by the random pairs $(F_k,M_k)$, $k=0, \ldots,n$. For every $n \ge 0$ we assume that, given $\mathcal F_n$, the pairs $(\rho_{ni}, \tau_{ni})$, $i\ge 1$,   are iid random variables with values in $\mathbb N_0\times \mathbb N_0$. In former investigations it has been assumed that their conditional distribution $\mu_n$ is non-random, here we allow that   $\mu_n$ depends on $L(F_n,M_n)$. Thus we deal with a population size dependent bisexual Galton-Watson process. Note that the random variables $X_n=L(F_n,M_n)$, $n =0,1,\ldots$, form a Markov chain with values in $\mathbb N_0$ and an absorbing state 0.

For the  bisexual Galton-Watson without population size dependence the question of recurrence/transience (or in other words the question, whether extinction appears with probability 1) has been completely solved for the large class of superadditive mating functions, see \cite{da} and the literature cited therein. In \cite{mo,mo2} the authors treated the case of mating functions depending on the population size. 

Here we consider the situation where the distribution of $(\rho_{ni},\tau_{ni})$ may depend on the number of couples $X_n$ (this is a particular case of the model introduced in \cite{mo3}). Then it is necessary to specify the function $L$ in more detail. We consider the prominent case
\[   L(x,y)= \min(x,ry)  ,\]
where $r\ge 1$ is a natural number ($r=1$ means monogamous mating and $r \ge 2$ polygamous mating). We restrict ourselves to the balanced situation when
\[ \mathbf E[\rho_{ni} \mid \mathcal F_n] = r\, \mathbf E[ \tau_{ni}\mid \mathcal F_n] \ \text{ a.s.} \]
 (the unbalanced case can be treated equally). In the case $r=1$ this means that in mean all females and males will find together in couples and it is only  due to random fluctuations that some will not succeed.

We like to apply our theorems to the process $(X_n)_{n \ge 0}$. The function $g$ evaluated on $x>0$  is given by \begin{align*} x+g(x)&= \mathbf E[ X_{n+1}\mid X_n=x  ] = \mathbf E_x\big[ \min\big(\sum_{i=1}^{x} \rho_{0i} , r \sum_{i=1}^{x}\tau_{0i}\big)  \big]  \ ,
\end{align*}
where we now use the notation $\mathbf E_x[\ \cdot\ ] =\mathbf  E[\ \cdot \mid X_0=x]$ familiar for Markov chains. Denoting
\[ e(x) = \mathbf E_x[ \rho_{01}  ]   \]
and using the identity $2\min(u,v)= u+v-|u-v|$, we may rewrite the above equation as
\[ g(x)+x = xe(x)  -\tfrac 12  \mathbf E_x[\big| \sum_{i=1}^{x}( \rho_{0i}-r\tau_{0i})\big| ] \] 
in the balanced case. The right-hand expectation can be asymptotically evaluated as follows: Assuming that the function
\[ v(x) = \mathbf{E}_x[ (\rho_{01}-r\tau_{01})^2  ] \]
has a finite, strictly positive limit for $x \to \infty$, i.e.
\[ v(x) \to \alpha >0 \quad \text{as } x \to \infty \ , \]
and assuming also
\begin{align} \mathbf{E}_x[ \rho_{01}^{2+\eta}+\tau_{01}^{2+\eta} ] \le c   \label{moment}
\end{align}
for some $\eta>0$, $c  <\infty$, we deduce from Lyapunov's version of the central limit theorem that
\[  \mathbf E_x[\big| \sum_{i=1}^{x}( \rho_{0i}-r\tau_{0i})\big|  ] =  \sqrt{ \alpha x}    ( \mathbf E[|N|] +o(1))   \quad \text{as } x \to \infty \ , \]
where $N$ has  a standard normal distribution. Since the right-hand expectation is equal to $\sqrt{2/\pi}$, we end up with
\[ g(x)=  (e(x)-1)x - \sqrt{\frac {\alpha x} { 2 \pi } }  + o(\sqrt x)  \quad \text{as } x \to \infty \ . \]

The  moments  of
$\xi_n$ can be obtained  as follows:
\begin{align*}
\mathbf E_x[ &|\xi_0|^{2+\delta}  ] \\ &= \mathbf E_x[ \big| \min\big(\sum_{i=1}^{x} \rho_{0i} , r \sum_{i=1}^{x}\tau_{0i}\big) - x-g(x) \big|^{2+\delta}  ]\\
&= \mathbf E_x[ \big| \min\big(\sum_{i=1}^{x} (\rho_{0i}-e(x)) ,  \sum_{i=1}^{x}(r\tau_{0i}-e(x))\big) + \sqrt{\frac {\alpha x} { 2 \pi } }  + o(\sqrt x)  \big|^{2+\delta}]
\end{align*}
From the Marcinkiewicz-Zygmund and the H\"older inequality we have
\begin{align*} \mathbf E_x[\big| \sum_{i=1}^x&(\rho_{0i}-e(x))\big|^{2+\eta}] \le c \mathbf E[\big(\sum_{i=1}^x |\rho_{0i}-e(x)|^2\big)^{1+\eta/2}] 
\\ &\le c \mathbf E[ \sum_{i=1}^x |\rho_{0i}-e(x)|^{2+\eta} x^{\eta/2}] = c x^{1+\eta/2} \mathbf E[|\rho_{01}-e(x)|^{2+\eta}]
\end{align*}
for some $c>0$. Similary the moment of the other sum may be estimated from above. Thus because of \eqref{moment} we may apply the multivariate central limit theorem to obtain for $0\le \delta < \eta$
\[ \mathbf E_x[ |\xi_0|^{2+\delta}] = x^{1+\delta/2} \mathbf E[ \big|\min(N_1,N_2)+ 
\sqrt{   \alpha/2 \pi }\big|^{2+\delta} ] +o(x^{1+\delta/2})   
\]
as $x \to \infty$, where the distribution of $(N_1,N_2)$ is bivariate normal. In parti\-cular the right-hand expectation is strictly positive. Since we are dealing with a time homogeneous Markov chain, \eqref{A2} is fulfilled for all $\delta < \eta$.

We are now ready to apply our theorems. The above formulas suggest to chose
\[ e(x)= 1+ \frac \beta {\sqrt x} + o(x^{-1/2}) \ \text{ as } x \to \infty  \] 
 for some real number $\beta$. Then asymptotically 
\[ g(x) \sim (\beta- \sqrt{\alpha/2\pi})\sqrt x \ , \quad \sigma^2(x)= \mathbf E_x [\xi_0^2] \sim \gamma x \]
for some $\gamma >0$. Thus applying Theorem 1 and 2 and noting that 0 is the only absorbing state, we end up with the following result.

\paragraph{Corollary.} {\em Assume that 
\[ e(x)= 1+ \frac \beta {\sqrt x} + o(x^{-1/2}) \ \text{ as } x \to \infty  \] 
and that $ \mathbf{E}_x[ \rho_{01}^{2+\eta}+\tau_{01}^{2+\eta} ] \le c$ for all $x$ large enough and some $\eta, c>0$.
 Then we have:
 \begin{enumerate}
 \item[{\em (i)}]
 If $\beta < \sqrt{ \alpha/2\pi}$, then the process $(X_n)_{n \ge 0}$ gets extinct with probability 1.
 \item[{\em (ii)}] 
 If  $\beta > \sqrt{ \alpha/2\pi}$, then $(X_n)_{n \ge 0}$ diverges  with positive probability.
 \end{enumerate}  }

\section{The general case}

Now we come back to the $d$-dimensional process $(X_n)_{n\ge0}$ with values in $\mathbb R_+^d$ satisfying
\[ X_{n+1}=MX_n + g(X_n)+ \xi_n \ \text{ with } \mathbf E[ \xi_n \mid \mathcal F_n]=0 \ , \ \mathbf E[ (\ell \xi_n)^2 \mid \mathcal F_n]= \sigma^2(X_n) \text{ a.s.}  \]
We recall that $M$ is assumed to be a primitive matrix, that is its entries are non-negative and there is a natural number $k$ such that the entries of $M^k$ are all strictly positive. We further assume that the Perron-Frobenius eigenvalue of $M$ is equal to 1. As is well-known, see \cite{se},  it has unique left and right eigenvectors  $\ell=(\ell_1,\ldots, \ell_d)$ and $r=(r_1, \ldots , r_d)^T$ with strictly positive components and normalized by $\ell_1+ \cdots + \ell_d= \ell_1r_1+ \cdots + \ell_dr_d=1$.

We like to obtain generalisations of the Theorems 1 and 2 above. First results in this direction are due to Adam \cite{ad} who derives recurrence and transience criteria in the case where $\ell g(x)$ and $\sigma^2(x)$ asymptotically behave like certain powers of $\ell x$.

 A first trial might be to look for suitable $d$-dimensional versions of the recurrence condition \eqref{recurr} resp. of the transience condition \eqref{trans} holding everywhere in the  $\mathbb R_+^d$ (up to a neighbourhood of the origin). However there occur interesting instances where in some parts of the state space one then would come across the recurrence condition and in others across the transience condition. 

Therefore we follow a different idea. From Perron-Frobenius theory it is known that the modulus of the other eigenvalues of $M$ are all smaller than 1. Roughly speaking this implies that after multiplication with $M$ vectors $x \in \mathbb R_+^d$  are pushes  towards the direction of the eigenvector $r$. The same effect is also active for the process $(X_n)_{n \ge 0}$ (provided it is not nullified by the presence of $g$ or the $\xi_n$). Thus on the event $\|X_n\|\to \infty$ one would expect that the sequence $X_n$ diverges asymptotically in the direction determined by $r$. Consequently the mentioned recurrence or transience conditions have to be required only in some vicinity of the ray $ \overline r=\{ \lambda r : \lambda \ge 0\}$ spanned by the vector $r$. This is the type of condition we are aiming at.

To make this intuition precise let us introduce some notation. Each vector $x \in \mathbb R^d$ can be uniquely dissected into two parts
\[ x= \hat x+ \check x\]
such that 
\[ \hat x \in \{\lambda r : \lambda \in \mathbb R\} \ \text{ and } \ \ell \check x= 0 \ . \]
From $\hat x= \lambda r$ it follows $\ell x= \lambda \ell r=\lambda $ because of $\ell r=1$. Thus
\[  \hat x = r\ell x\ . \]
Now we may characterize the vectors $x$ close to the ray $\overline r$ by the property that $\| \check x \|$ is small compared to $\|x\|$.

Similary the process $(X_n)$ can be splitted:
\[X_n = \hat X_n+ \check X_n \ . \]

We are now ready to formulate our results. The assumptions are analoguous to the 1-dimensional case. The condition of near-criticality   now reads
\begin{align}   \|g(x)\| = o(\|x\|) \quad \text{as } \|x\| \to \infty \tag{A1$^*$} \label{A11} 
\end{align}
and the condition of moment boundedness gets the following form: There are $c>0$ and $\delta >0$ such that for $p=2+\delta$ and  $\| X_n\| \ge c$
\begin{align}   \mathbf E[ \|\xi_n\|^{p} \mid \mathcal F_n] \le c \sigma^{p}(X_n) \ .
\tag{A2$^*$} \label{A22} 
\end{align}

For convenience we formulate here our theorems only for the case that $g(x)$ has only non-negative components, otherwise the required conditions are somewhat more involved.

\paragraph{Theorem 3.} {\em  Let \eqref{A11} and \eqref{A22} be fulfilled and let $g(x) \ge 0$ (component\-wise) for $\|x\|$ sufficiently large. Assume that there is an $\varepsilon >0$ and that for any $a>0$ there is some $b>0$ such that for all $x \in \mathbb R_+^d$ we have
\begin{align} \|x\| \ge b \ , \|\check x\|^2 \le a \|x\|\cdot \|g(x)\| \quad \Rightarrow \quad \ell x \cdot \ell g(x) \le \frac{1-\varepsilon}2 \sigma^2(x) \ . \label{cond}
\end{align}
Then
\[ \mathbf P( \|X_n\| \to \infty \text{ as } n \to \infty)= 0 \ . \]}

Note that in view of \eqref{A11} the requirement $\|\check x\|^2 \le a \|x\|\cdot \|g(x)\|$ indeed defines a region within $\mathbb R_+^d$ which is located in the vicinity of the ray $\overline r$ (since $x \in \overline r$ implies $\check x=0$). Note also that in the 1-dimensional case the requirement \eqref{cond} reduces to \eqref{recurr}, because then $\check x=0$ and $\ell x =x $ for any $x \in \mathbb R_+$.

\paragraph{Theorem 4.} {\em  Let \eqref{A11} and \eqref{A22} be fulfilled and let $g(x) \ge 0$  for $\|x\|$ sufficiently large. Also let $\sigma^2(x)$ be bounded away from zero for all $x \in \mathbb R_+^d$ with $u< \ell x < v$, where $0<u<v<\infty$, and let
\begin{align}  \sigma^2(x)= O(\|x\|^2 \log^{-2/\delta} \|x\|) \ \text{ for } x \to \infty \text{ and for some } \delta >0 \ . \label{sigma2} \end{align}
Assume that there is an $\varepsilon >0$ and that for any $a>0$ there is some $b>0$ such that for all $x \in \mathbb R_+^d$ we have
\begin{align*} \|x\| \ge b \ , \|\check x\|^2 \le a \sigma^2(x) \quad \Rightarrow \quad \ell x \cdot \ell g(x) \ge \frac{1+\varepsilon}2 \sigma^2(x) \ .  
\end{align*}
Then there is a number $m< \infty$ such that
\[ \mathbf P( \limsup_n \|X_n\| \le m \text{ or } \lim_n \|X_n\| = \infty)=1 \ . \]
If additionally for every constant $c>0$ there is a natural number $n$ such that we have $\mathbf P(\|X_n\|>c, c\|\check X_n\|\le \|X_n\|) >0$, then
\[ \mathbf P(\|X_n\| \to \infty \text{ for } n \to \infty)>0\] and \[\mathbf P(\|\check X_n\|= o(\|X_n\|) \mid \|X_n\|\to \infty )=1 .\]}

Now it is the condition $\|\check x\|^2 \le a \sigma^2(x)$ which in view of \eqref{sigma2} defines a vicinity of the ray $\overline r$. The last statement says that the process diverges in the direction of the ray $\overline r$. Again for $d=1$ this Theorem reduces completely to Theorem 2. 

Note that we did not specify which norm $\| \cdot \|$ we used. The choice makes no difference because as is well-known all norms on $\mathbb R^d$ are equivalent in the sense that for two norms $\| \cdot  \|_1$ and $\|  \cdot \|_2$ there are number $c_1,c_2>0$ such that $c_1\| \cdot  \|_1 \le \| \cdot  \|_2 \le c_2 \|  \cdot \|_1$. Thus the formulated conditions and statements do not depend on the choice of the norm.

The proof of both theorems will be given elsewhere.

\section{Example: The multivariate GW-process}

In the multivariate  Galton Watson process  $(X_n)_{n \ge 0}$ each generation consists of $d$ different types of individuals, thus $X_n= (X_{n,1}, \ldots, X_{n,d})^T$. Then 
\[ X_{ n+1 }= \sum_{j=1}^d \sum_{i=1}^{X_{n,j}} \zeta_{nij }  \ \text{ or }\ X_{ n+1,k}= \sum_{j=1}^d \sum_{i=1}^{X_{n,j}} \zeta_{nijk}  \]
with $\zeta_{nij}=(\zeta_{nij1}, \ldots, \zeta_{nijd})^T$. Here $\zeta_{nijk}$ is considered to be the   offspring number of individuals of type $k$ born by the $i$-th individual of type $j$ in generation $n$. It is assumed that for $n\ge 0$ and  given $X_0, \ldots,X_n$ the random vectors $\zeta_{nij}$, $i,j \ge 1$, are independent with a distribution which may depend on $j$ but neither  on $i$ nor on $n$. We allow that this distribution depends also on $X_n$, then $(X_n)_{n \ge 0}$ is a Markov chain  and constitutes a population size dependent multivariate Galton-Watson process as introduced by Klebaner, see e.g. \cite{kleb}. For a model where $M$ is no longer primitive, compare Jagers and Sagitov \cite{ja}.

Now given the state $x=(x_1, \ldots, x_d)^T \in \mathbb R_+^d$ expectations are determined as 
$ \mathbf E_x[ X_{1 } ] = \sum_{j=1}^d \mathbf E_x[\zeta_{01j }] x_j $ or
\[ \mathbf E_x[X_1]  = E_x x \] with the $d\times d$ matrix of expectations \[ E_x = (\mathbf E_x[\zeta_{011}], \ldots, \mathbf E_x[ \zeta_{01d}])\ .\]
In our context this means that
\[ g(x)= (E_x-M)x \]
with some primitive matrix $M$ and
\[ \xi_n =  \sum_{j=1}^d \sum_{i=1}^{X_{n,j}} (\zeta_{nij}- \mathbf E_x[\zeta_{01j }] )\ . \]
By conditional independence we get
\begin{align*}
\sigma^2(x)= \mathbf E_x[ (\ell \xi_0)^2 ] = \sum_{j=1}^d  \mathbf{Var}_x(\ell \zeta_{01j})x_j
\end{align*}
or
\[ \sigma^2(x)= \ell \Gamma_x \ell^T \]
with
\[ \Gamma_x = \sum_{j=1}^d \mathbf{Cov}_x (\zeta_{01j})x_j \]
and the $d\times d$ covariance matrices $\mathbf{Cov}_x(\zeta_{01j})$  of $\zeta_{01j}$, $j=1,\ldots,d$.

As to condition \eqref{A22} it is suitable here to work with the $\ell_1$-norm $\|x\| = \sum_k |x_k|$.
Then for $p=2+\delta>2$
\begin{align*}
\| \xi_n\|^p =  \big(\sum_{k=1}^d|\sum_{j=1}^d  \sum_{i=1}^{X_{n,j}} (\zeta_{nijk}- \mathbf E_x[\zeta_{01jk }] |\big)^p\le d^{2p} \sum_{j,k=1}^d \big| \sum_{i=1}^{X_{n,j}} (\zeta_{nijk}- \mathbf E_x[\zeta_{01jk })] \big|^p
\end{align*}
Applying again the Marcinkiewicz-Zygmund and the H\"older inequality we obtain, for some $c>0$,
\begin{align*}
\mathbf E_x[\|\xi_0\|^p] \le c d^{2p} \sum_{j,k=1}^d x_j^{p/2} \mathbf E_x [|\zeta_{01jk}- \mathbf E_x[\zeta_{01jk }]|^p]\ .
\end{align*}
Assuming now that there is a number $b>0$ such that for all $x \in \mathbb R_+^d$
\begin{align} \mathbf E_x [|\zeta_{01jk}- \mathbf E_x[\zeta_{01jk }]|^p] \le b  \mathbf E_x [(\zeta_{01jk}- \mathbf E_x[\zeta_{01jk }])^2]^{p/2} \label{condition}
\end{align}
we obtain
\begin{align*}
\mathbf E_x[\|\xi_0\|^p]  \le bcd^{2p} \big( \sum_{j,k=1}^d x_j \mathbf E_x [(\zeta_{01jk}- \mathbf E_x[\zeta_{01jk }])^2]\big)^{p/2} \ ,
\end{align*}
that is
\[ \mathbf E_x[\|\xi_0\|^p]  \le bcd^{2p} (\text{trace}\ \Gamma_x)^{p/2} \ . \]
Thus to attain validity of \eqref{A22}  we require besides \eqref{condition} that there is a constant $c>0$ such that for all $x \in \mathbb R^d_+$ we have
\[ \text{trace } \Gamma_x \le c\, \ell \Gamma_x \ell^T \ . \] 
This is fulfilled if e.g. the covariance matrices $\mathbf{Cov}_x (\zeta_{01j})$ have only non-negative entries, but it may fail in general.

Now we are ready to apply our theorems to special cases as those discussed by Klebaner \cite{kleb} and Adam \cite{ad}. Details are left to the reader.


\end{document}